%Z input Amstex
\input amstex
\documentstyle{amsppt}
  \catcode`@11
%\outer\def\head#1\endhead{\relaxnext@ \def \nextii@ {\bigbreak 
%\bgroup \let \\=\cr \global \setbox \headingbox@ \vbox \bgroup 
%\tabskip \z@ \filhss@ \halign to\hsize \bgroup \tenpoint \smc 
%\hfil \ignorespaces ####\unskip \hfil \cr }\overlong@ \futurelet 
%\next \next@  #1 \endheading}

%\let\remark\proclaim
%\let\endremark\endproclaim
%\let\definition\proclaim
%\let\enddefinition\endproclaim

%%%%%  some stuff that was commented out in amstex.tex:
%%%  not needed anymore Aug 26, 91
%  End definition of Euler Fraktur font.

\catcode`@12

\let\olddate=\date
\def\date#1{\olddate#1\enddate}
\let\oldthanks\thanks
\def\thanks#1{\oldthanks#1\endthanks}
\let\oldabstract\abstract
\def\abstract#1{\oldabstract#1\endabstract}

\def\martinstricks{}   

\newtoks\fooline
\def\today{\number\day.\number\month.\number\year}
\fooline={{\tt\jobname}\hfill{\tenrm\relax\folio\hfill\today}}

\catcode`@11

\newtoks\oldoutput
\oldoutput=\the\output
\output={\setbox255\vbox to \ht255 
         {\unvbox255 %
\iffirstpage@\else 
\smallskip\line{\the\fooline}\vss\fi}\the\oldoutput}

\catcode`@12

\martinstricks

\def\demo#1{\smallskip\noindent{\smc #1}}
\def\enddemo{\medskip}

\hsize=4.5 in
\baselineskip=18pt
\magnification=1200
\topmatter
\title Existence of Endo-Rigid Boolean Algebras
\endtitle
\endtopmatter
\bigskip
\bigskip
\document
In [Sh 2] we, answering a question of Monk, have
explicated the notion of ``a Boolean algebra with no
endomorphisms except the ones induced by ultrafilters on
it'' (see \S 2 here) and proved the existence of one with
character density $\aleph_0$, assuming first
$\diamondsuit_{\aleph_1}$ and then only $CH$.  The idea
was that if $h$ is an endomorphism of $B$,not among the
``trivial'' ones, then there are pairwise disjoint
$D_n\in B$ with $h(d_n)\not\subset d_n$.  Then we can,
for some $S\subset \omega$, add an element $x$ such that
$d_n\in B$ with $h(d_n)\not\subset d_n$>  Then we can,
for some $S\subset\omega$, add an element $x$ such that
$d\leq x$ for $n\in S, x\cap d_n=0$ for $n\not\in S$ while
forbidding a solution for $\{y\cap h(d_n):n\in S\}\cup
\{y\cap h(d_n)=0:n\not\in S\}$.  Further  analysis showed
that the point is that we are omitting positive
quantifier free types.  Continuing this Monk succeeded to
prove in $ZFC$, the existence of such Boolean algebras of
cardinality $2^{\aleph_0}$.  In his proof he 

\roster\item"(a)" replaced some uses of the countable
density character by the $\aleph_1$-chain condition
\item"(b)" generally it is hard to omit $<2^{\aleph_0}$
many types but because of the special character of the
types and models involve, using $2^{\aleph_0}$ almost
disjoint subsets of $\omega$, he succeeded in doing this
\item"(c)"  for another step in the proof (ensuring
indecomposability - see Definition 2.1) he ( and
independently by Nikos) found it is in fact easier to do
this when for every countable $I\subset B$ there is $x\in
B$ free over it.
\endroster

The question of the existence of such Boolean algebras in
other cardinalities open (See [DMR] and a preliminary
list of problems for the handbook of Boolean Algebras by
Monk).

We shall prove (in $ZFC$) the existence of such $B$ of
density character $\lambda$ and cardinality
$\lambda^{\aleph_0}$ whenever $\lambda >\aleph_0$.  We can
conclude answers to some other questions from Monk's
list, (combine 3.1 with 2.5).  We use a combinatorial
method from [Sh 3],[Sh 4], that is represented in section
1.

In [Sh 1], [Sh 6] (and [Sh 7]) the author offers the
opinion that the combinatorial proofs of [Sh 1], Ch VIII
(applied thee for general first order theories) should be
useful for proving the existence of many non-isomorphic,
ad/or pairwise non-embeddable structures which have few
(or no) automorphisms or endomorphism of direct
decompositions etc.  As an illumination, in [Sh 6] a
rigid Boolean algebra in every $\lambda^{\aleph_0}$
omitting countable types along the way, the method is
proved in $ZFC$, nevertheless it has features of the
diamond.  It has been used  (so in Gobel and Corner [CG]
and Gobel and Shelah [GS1],[GS2].  See more on the method
and on refinements of it in [Sh 4] and [Sh 3] and mainly
[Sh 5].
\bigskip
\heading The combinatorial pronciple\endheading
\bigskip
{\bf Content}  Let $\lambda>k$ be fixed
infinite cardinal 
\bigskip
We shall deal with the case $cf\, \lambda\aleph_0,
\lambda^{\aleph_0} =\lambda^k$, and usually
$k=\aleph_0$.  Let $\Cal M$ be a function symbols, each
with $\leq k$ places, of power $\leq\lambda$. Let
$\zeta(\lambda,k)$ be the vocabulary with function
symbols $\{Fi,j:i<\lambda \ j<k\}$ where $F_i,j$ is a
$j$-place function symbol.  Let
$\Cal M=\Cal M_{\lambda,k}(T)$
the free $\zeta(\lambda,k)$-algebra generated by
$T\buildrel\text{def}\over =
{}^{k>}\lambda(=\{\eta:\eta$ a sequence $<k$ of
ordinals $<\lambda$) (We could have as well considered
$T$ as a set of urelements, and let $\Cal M$ be the family
$H_{<k'} (T)$ of sets hereditarily of cardinality
$<k$ build from urelements].

\bigskip
\proclaim{1.2 Explanation}  We shall let $B_0$ be the
Boolean Algebra freely generated by $\{\eta:\eta\in T\},
B^C_0$ its completion and we can interprete $B^c_0$ as a
subset og $\Cal M$ (each $a\in B^c_0$ has the form
$\cup_{\eta<\omega} \tau_\eta$ where $\tau_\eta$ is a
Boolean combination of members of $T$, so as we have in
$L\, \aleph_0$-place function symbols there is no
problem).  As the $\eta\in T$ may be over-used we
replaced them for this purpose by ?? (e.g. let $F\in L$
be a monadic function symbol, $x_\eta=F(\eta)$).

Our desired Boolean Algebra $B$ will be a subalgebra of
$B^c_0$ containing $B_0$.\endproclaim
\bigskip
\definition{1.3 Definition}
\roster
\item Let $L_n$, for $n<w$, be fixed vocabularies
(= signatures), $|L_n|\leq L_n k, L_{n+1} \subset
L_{n+1}$, (with each predicate function symbol finitary
for simplicity) let $P_n\in L_{n+1}-L_n$ be monadic
predicates.
\item  Let $\Cal F_n$ be the family of sets (or
sequences) of the form $\{(f_\ell,N_\ell):\ell\leq\}$
satisfying

\indent a) $f_\ell:{}^{\ell\geq}k\rightarrow T$ is
a {\it tree emdedding} i.

\ \ \ \ \ (i)  $f_\ell$ is length preserving
i.e.$\eta,f_ell(\eta)$ have the same length.

\ \ \ \ \ (ii)  $f_\ell$ is order preserving i.e. for
$\eta,\nu\in {}^{\ell\leq k},\eta <\nu$ iff 
$f_\ell(\eta)<f_\ell(\nu)$.

\indent b) $f_{e+1}$ extend $f_\ell$ (when $\ell+1\leq n$)

\indent c) $N_\ell$ is an $L'_\ell$-model of power $\leq
k,|N_\ell|\subset |XX|$,where $L'_\ell\subset
L_\ell$.

\indent d) $L'{\ell+1} \cap L_\ell=L'_\ell$ and
$N_{\ell+1}\restriction L'\ell$ extends $N_\ell$

\indent e) if $P_m\in L'_{m+1}$, then
$P^{N_\ell}_m=|N_m|$ when $m<\ell\leq n$ and

\indent f) $Rang(f_\ell) -\cup_{m<\ell} Rang(f_m)$ is
included in $N_\ell|-\cup_{m<\ell}|N_m|$.

\item  Let $\Cal F_w$ be the family of pairs $(f,N)$ such
that for some $(f_\ell,N_\ell)(\ell< \omega)$ the following
holds:

\indent (i) $\{(f_\ell,N_\ell):\ell\leq n\}$ belongs to
$\Cal F_n$ for $n<w$.

\indent (ii) $f=\cup_{\ell<\omega} f_\ell,N=\cup_{n<w}
N_n$, (i.e. $|N|=\cup_{n<\omega}|N_n|, L(N)=\cup_n L(N_n)$,
and $N\restriction L(N_n)=\cup_{n<m<\omega}
N_m\restriction L(N_m)$
\item For any $(f,N)\in \Cal F_\omega$ let $(f_n, N_n)$ be
as above (it is easy to show that $(f_n,N_n)$ is uniquely
determined - notice d),e) in (2),) so for
$(f^\alpha,N^\alpha)$ we get $(f^\alpha_n,N^\alpha_n)$
\item A branch of $Rang(f)$ or of $f$ (for $f$ as in
(3)) is just $\eta\in {}^\omega\lambda$ such that for every
$n<\omega,\eta \restriction n \in Rang(f)$.
\endroster
\enddefinition
\bigskip
{\bf 1.4 Explanation of our Intended Plan}
 (of Constructing e.g. the Boolean algebra)
\bigskip
We will be given
$W=\{(f^\alpha,N^\alpha):\alpha<\alpha(^*)\}$, so that
every branch $\eta $ of $f^\alpha$ converges to some
$\zeta(\alpha),\zeta(\alpha)$ non decreasing (in
$\alpha$).  We have a free object generated by $T{(B_0}$ in
our case) and by induction on $\alpha$ we define
$B_\alpha$ and $a_\alpha, B_\alpha$ increasing
continuos, such that $B_{\alpha+1}$ is an
extension of $B_\alpha,a_\alpha\in B_{\alpha+1}-B_\alpha$
usually $B_{\alpha+1}$ is generated by $B_\alpha$ and
$a_\alpha$ is in the completion of $B_0$).   Every
element will depend on few $(\leq k)$ members of
$T$, and $a_\alpha$ ``depends'' in a peculiar way: the
set $Y_\alpha \subseteq T$ on which it ``depends'' is
$Y^0_\alpha \cup Y^1_\alpha$ where $Y^0_\alpha$ is
bounded below $\zeta(\alpha)$ (i.e. $Y^0_\alpha\subseteq
{}^{\omega>}\zeta$ for some $\zeta<\zeta(\alpha))$ and
$Y^1_\alpha$ is a branch of $f^\alpha$ or something
similar.  See more in 1.8.
\bigskip
\definition{1.5 Definition of the Game}  We define for
$W\subseteq \Cal F_\omega$ a game $Gm(W)$, which asts
$\omega$-moves.

In the $n$-th move:

{\it Player II}: \   Choose $f_n$, a tree-embedding of
${}^{n\geq}\lambda$, extending $\cup_{\ell<n} f_\ell$,
such that $Rang(f_n)-\cup_{\ell<n}\, Rang(f_\ell)$ is
disjoint to $\cup_{\ell<n}|N_\ell|$; then

{\it player I} chooses $N_n$ such that
$\{(f_\ell,N_\ell):\ell\leq n\}\in \Cal F_n$.  

In the end player II wins if
$(\cup_{n<\omega}f_n,\cup_{n<\omega}N_n)\in W$.

\enddefinition
\bigskip
\remark{1.6 Remark}  We shall be interested in $W$ such
that player II wins (or at least does not lose) the game,
but $W$ is ``thin''.  Sometimes we need a strengthening
of the first player in two respects: he can force (in the
$n$-th move) $Rang(f_{n+1})-Rang(f_n)$ to be outside a
``small'' set, and in the zero move he can determine an
arbitrary initial segment of the play.
\endremark
\bigskip
\definition{1.7 Definition}  We define, for $W\subseteq
\goth F_\omega$, a game $GM'(W)$ which lasts
$\omega$-moves.

In the zero move

{\it player II} choose $f_0$, a tree embedding of
${}^{0\geq}$ to ${}^{0\geq}\lambda$ (but there is only
one choice).
\smallskip
{\it player I} chooses $k<\omega$ and
$\{(f_\ell,N_\ell):\ell\leq \}\in \goth f_k$, and
$X_0\subset T,|X_0|<\lambda$.  In the $n$-th move, $n<0$:
\smallskip
{\it player II} chooses $f_{k+n}$ a tree embedding
of ${}^{(k+n)\geq}k$ into
${}^{(k+n)\geq}\lambda$, with $Range\,
f_{k+n}-\cup_{\ell<k+n} Rang\, f_\ell$ disjoint to
$\cup_{\ell<k+n} N_\ell\cup \cup_{\ell<n} X_\ell$
\smallskip
{\it player I} choose $N_{k+n}$ such that
$\{(f_\ell,N_\ell):\ell \leq k+n\}\in \goth
F_{k+n}$ and $X_n\subseteq T,|X_n|<\lambda$
\enddefinition
\bigskip
\remark{1.8 Remark}  What do we want from $W$?: First
that by adding an element (to $B_0$) for each $(f,N)$ we
can ``kill'' every undesirable endomorphism, for this it
has to encounter every possible endomorphism, and this
will be served by ``$W$ a barrier''.  For this $W=\Cal
F_\omega$ is O.K. but we also want $W$ to be thin enough so
that various demands will have small interaction, for
this disjointness and more are demanded.
\endremark
\bigskip
\definition{1.6 Definition} 
\roster \item We call $W\subseteq \Cal F_\omega$ a {\it
strong barrier} if player II wins in $gm(W)$ and even
$Gm'(W)$ (which just means he has a winning strategy.)
\item  We call $W$ a {\it barrier} if player    I does
not win in $Gm(W)$ and even does not win in $Gm'(W)$.
\item  We call $W$ disjoint {\it if} for any distinct
$(f^\ell,N^\ell)\in W(\ell=1,2), f^1$ and $f^2$ have no
common branch.
\endroster
\enddefinition
\bigskip
\proclaim{1.7 The Existence Theorem}
\roster
\item if $\lambda^{\aleph_0}=\lambda^k, cf
\lambda>\aleph_0$ then there is a strong disjoint barrier.
\item Suppose $\lambda^{\aleph_0}=\lambda^k,\, cd \,
\lambda>\aleph_0$.  Then there is
$W=\{(f^\alpha,N^\alpha):\alpha<\alpha^*\}\subseteq \Cal
F_w$ and a function $\zeta:\alpha^*\rightarrow \lambda$
such that:
\smallskip
\ \ \ \ (a)  $W$ is a strong disjoint barrier, moreover
for every stationary $S\subset
\{\delta<\lambda:cf\,\delta=\aleph_0\}\, \{(f^\alpha,
N^\alpha):\alpha<\alpha^*,\zeta (\alpha)\in S\}$ is a
disjoint barrier.
\ \ \ \ (b)  $cf\,(\zeta(\alpha))=\aleph_0$ for
$\alpha<\alpha^*$.
\ \ \ \ (c)  Every branch of $f^\alpha$ is an increasing
sequence converging to $\zeta(\alpha)$.
\ \ \ \ (d)  If $\overline \eta$ is a sequence from $T$ (of
any length $\gamma<k^*), \tau(\bar x)$ a term of
$L$and every $X_i$ really appears, $lg(\overline
x)=\gamma$ and $\tau(\overline \eta) \in N^\alpha$ then
$\bar \eta\subseteq N^\alpha \cap T$. 
\ \ \ \ (e)  If
$\zeta(\beta)=\zeta(\alpha),\beta+k^{\aleph_0}\leq
\alpha<\alpha^*$ and $\eta$ is a branch of $f^\alpha$ then
$\eta\restriction k\not\in N^\beta$ or some
$k<w$.
\ \ \ \ (f)  If $\lambda=\lambda^k$ we can demand: 
if $\eta$ is a branch of $f^\alpha$ and $\eta \restriction
k\in N^\beta$ for all $k<\omega$ (where
$\alpha,\beta<\alpha^*$) then $N^\alpha\subseteq N^\beta$
(and even ??? $N^\alpha_n\in N^\beta$ if
$\Cal M=H_{<k^+}(T))$.
\endroster
\endproclaim
\bigskip
\bigskip
\heading Reference \endheading
\bigskip
{\bf 2 Preliminaries on Boolean Algebras}
\smallskip
We review here some easy material from [Sh 2].

\definition{2.1 Definition} 
\roster
\item  For any  endomorphism $h$ of a Boolean Algebra
$B$. let $EX \, Ker(h)=\{x_1\cup x_2:\,h(x_1)=0$, and
$h(y)=y$ for every $y\leq x_2\}$.

$Ex\, Ker^*(h)=\{x\in B:$ in $B/Ex\, Ker(h)$, below
$x/Ex\, Ker^*(h)$, there are only finitely many
elements$\}$.
\item  A Boolean algebra is endo-rigid {\it if} for every
endomorphism $h$ of $B,B/Ex\, Ker(h)$ is finite
(equivalently: $1_B\in Ex\, Ker(h))$.
\item  A Boolean algebra is indecompensable {\it if}
there are no two disjoint ideal  $?_0, I_1$ of $B$, each
with no maximal member which generate a maximal ideal
$\{a_0\cup a_1:\,a_0\in I_0,a_1\in I_1\}0$.
\item  A Boolean algebra $B$ is $\aleph_1$-compact if for
pairwise disjoint $?_n\in B(n<\omega)$ for some $x\in B,
x\cap d_{2n+1}=0, x\cap d_{2n}=d_{2n}$.
\endroster
\enddefinition
\bigskip
\proclaim{2.2 Lemma}
\roster
\item  A Boolean algebra $B$ is endo-rigid {\it iff} for
every endomorphism of $B$is the endomorphism of some
scheme (see Definition 2.3 below).
\item  A Boolean algebra $B$ is endo-rigid and
indecomposable {\it iff} every endomorphism of $B$ is the
endomorphism of some simple-scheme (see Def 2.3 below).
\endroster 
\endproclaim
\bigskip
\definition{2.3 Definition} 
\roster\item"(1)" A scheme of an endomorphism of $B$
consists of a partition $a_0,a_1,b_0,\dots,b_{n-1}$,

$c_0,\dots,c_{m-1}$ of $B$ of maximal nonprincipal ideal
$I_\ell$ below $b_\ell$ for $\ell<n$, nonprincipal
disjoint ideals $I^0_\ell,I^1_\ell$ below $c_\ell$ for
$\ell<m$, a number $k<n$, and a partition
$b^*_0,\dots,b^*_{n-1}, c^*_0,\dots,c^*_{m-1}$ of
$a_0\cup b_0\cup ,\dots\cup b_{k-1}$.  We assume also
that $[k+m>0\Rightarrow a_0=0],
[(n-k)+m>0\Rightarrow a_1=0]$ and except in those
cases there are no zero elements in the partition.
\item"(2)" The scheme is simple if $m-0$.
\item"(3)" The endomorphism of the scheme is the unique
endomorphism $T:B\rightarrow$ such that:
\indent (i) $Tz=0$ when $x<a_0$ or $x\in
I_\ell,\ell<k$, or $x\in I^0_\ell,\ell<m$.

\indent (ii) $Tz=x$ when $x\leq a_1$ or $x\in
I_\ell,k\leq \ell<n$ or $x\in I^1_\ell,\ell<m$.

\indent (iii) $T(b_\ell)=b_\ell^*$ when $\ell<k$.

\indent (iv) $T(b_\ell)=b_\ell\cup b^*_\ell$ when
$k\leq \ell<n$.

\indent (v)  $T(C_\ell)=c_\ell\cup c^*_\ell$ when
$\ell<m$.

\endroster
\enddefinition
\bigskip
\proclaim{2.4 Claim}  If $h$ is an endomorphism of a
Boolean Algebra $B$, and $B/ Ex\, Ker(h)$ is infinite
{\it then} there are pairwise disjoint $d_n\in
B(n>\omega)$ such that $h(d_n) \not\subset d_n$.  By easy
manipulation  we can assume that $h(d_n) \cap d_{n+1}
\not= 0$, and if $B$ satisfies the c.c.c then $\{d_n:\,
n<w\}$ is a maximal antichain.  So, an endomorphism of a
scheme is a ``trivial'' endomorphism generated by ideals.
\endproclaim
\bigskip
\proclaim{2.5 Lemma}  
\roster
\item  Every endo-rigid Boolean Algebra $B$ is a Hopfian
and dual Hopfian.  Even $B+B$ is Hopfian (and dual
Hopfian) but not rigid.
\endroster
\endproclaim
\bigskip
\demo{Proof}  Easy to check using 2.2, 2.3.
\enddemo
\bigskip
\bigskip
\heading The Construction\endheading
\bigskip
\proclaim{3.1 Main Theorem}  Suppose $\lambda>\aleph_0$. 
Then there is a B.A. (Boolean Algebra) $B$ such that:
\roster
\item
 $B$ satisfies the c.c.c.
\item $B$ has power $\lambda^{\aleph_0}$ and
$T_1(B)=\lambda$ where $T_1$ is the density character.
\item $B$ is endo-rigid and indecomposable.
\endroster
\endproclaim
\bigskip
\demo{Proof}  We concentrate on the case $cf(\lambda)\geq
\aleph_1$ (on the case $cf\, \lambda=\aleph_0$ see [Sh 5,
\S 2, \S 3]) we shall use Theorem 1.7, and let
$W=\{(f^\alpha,N^\alpha):\, \alpha<\alpha^*\}$, the
function $\zeta, \Cal M$ and $T={}^{\omega>}\lambda$ be as
there.

We will think of the game as follows:  player I tries to
produce a non trivial endomorphisms $h$.  Player II
supplies (via range  ($f_1$)) elements is $B-0$ and
challenges player I for defining $h$ on them.  So player
I plays models $N_i$ in the language $L_i$ with a
distinguished function symbol $h$ which is an
endomorphism of Boolean Algebras.  In the end, as $W$ is
a strong barrier, we will get a model $N^\alpha\in W$, in
the language $\cup_{i<\omega} L_i$ which includes a
function $h$.  We can think of $N^\alpha$ as a Boolean
Algebra $\subseteq B^c_0$ with an endomorphism $h$.
\enddemo
\bigskip
{\bf Stage A}  Let $B_0$ be the B.A. freely generated by
$\{x_\eta:\, \eta\in T\}$, let $x_\eta=a_\eta$ and
$B^c_0$ be its completion.  For $A\subseteq B^c_0$ let
$\langle A\rangle_{B^0_c}$ can be represented as a
countable union of members of $B_0$, so $w.o.l.g.\,
B^c_0\subseteq \Cal M$.  We say $x\in B^c_0$ is based on
$J\subseteq  {}^{>\lambda}$ if it is based on $\{x_\nu:\,
\nu\in J\}$ [i.e.$x=\cup_n y_n$, each $ y_n$ is
in the subalgebra generated by $\{x_\nu:\, \nu\in J\}$]
and let $\underline d(x)$ be the minimal such $J$.  We
shall now define by induction on $\alpha <\alpha^*$, the
truth values of ``$\alpha \in J$'',$\eta_\alpha$, and
members $\alpha_\alpha, b^\alpha_n, \,
c^\alpha_m,d^\alpha_m,\tau^\alpha_m$ of $B^c_0$ such
that, letting $B_\alpha=\langle
B_0,\alpha_i\,i<\alpha,i\in J\rangle_{B^c_0}$:

\roster 
\item $\eta_\alpha$ is a branch of
$Rang\,(f^\alpha),\eta_\alpha\not= \eta_\beta$ for
$\beta<\alpha$
\item if $\alpha\in J$, then for some $xi<\zeta(\alpha)$:

$a_\alpha=\cup_m(\tau^\alpha_m\cap d^\alpha_m)$ where
$\langle d^\alpha_m:\,m<\omega\rangle$ is a maximal
antichain of non zero elements (of $B^c_0$)
$\cup_m\underline d(d^\alpha_m)\subseteq
{}^{\omega>}\xi,\tau^\alpha_m\in \langle x_p:\,
\eta_\alpha\restriction m\leq  p,p\in T\rangle_{B^c_0}$,
and $\tau^\alpha_m\cap d^\alpha_m> 0$. 
\item  if
$\alpha\in J$, then $b^\alpha_n,d^\alpha_n\in
N^\alpha_0,\, c^\alpha_n,\tau^\alpha_m\in N^\alpha$ (hence
each is based on $\{x_\nu:\,\nu\in {}^{\omega>},\nu\in
N^\alpha\})$, and $b^\alpha_n\cap b^\alpha_m=0$ for
$n\not=m$.
 \item for $\beta<\alpha,\, \beta\in J,
B_\alpha$ omits $p_\beta=\{x\cap b^\beta_n=c_n^\beta:\,
\eta<\omega\}$. \endroster
\enddemo
\bigskip
\remark{Remark}  Many times we shall write
$\beta<\alpha<\alpha^*$ or $w\subseteq
\alpha<\alpha^*$ instead $\beta\in \alpha\cap
J,w\subseteq \alpha \cap J$.
\endremark
Before we carry the construction note:
\bigskip
\proclaim{3.2 Crucial Fact}:  For any $x\in B_\alpha$ there
are $k,\xi<\zeta$, and $\alpha_0<\dots<\alpha_k$
such that
$\zeta(\alpha_0)=\zeta(\alpha_1)=\zeta(\alpha_2)=\dots=\zeta(\alpha_k)=\zeta,x$
is based on $\{x_\nu:\, \nu\in {}^{\omega>}\xi$ {\it or}
$\nu\in \underline d(\tau^{{\alpha_\ell}_m}$, for some
$\ell\geq k,m<\omega\}$.
\endproclaim
\bigskip
{\bf Stage B}  Let us carry the construction.  For
$\xi<\lambda,w \subseteq \alpha^*$ let

$$I_{\xi,w}=\{\nu:\, \nu\in {}^{w>}\xi
\text{or}\,\nu\in\cup_{{m<w}_{\gamma\in w}}\underline
d(\tau^{\gamma_m})\}$$

We let $\alpha \in J$ {\it iff} $|N^\alpha|\subseteq
B_\alpha,N^\alpha=(B^c_0\restriction |N^\alpha|,h_\alpha)$
where $h_\alpha$ is an endomorphism of $B^c_0
\restriction |N^\alpha|$ hence maps $N^\alpha_n$ into
$N^\alpha_n$ for $n<\omega$) and there are $d^\alpha_m\in
N^\alpha_0$ for $m<\omega, d^\alpha_m\not=0,
d^\alpha_m\cap d^\alpha_\ell=0$ for $m\not=\ell$, such
that for some $xi<\zeta(\alpha)$  each $d^\alpha_m$ is
based on ${}^{\omega>}\xi$, and there are a branch
$\eta_\alpha$ of $Rang\,(f^\alpha)$ and $\tau^\alpha_m\in
N^\alpha(m<\omega)$ as in 1),2) above, such that if we add
$\cup_{n<\omega}(\tau^\alpha_n\cap d^\alpha_\ell)$ to
$B_\alpha$, each $p_\beta(\beta<\alpha$ is still omitted
as well as $P_\alpha=\{x\cap
h_\alpha(d^\alpha_m)=h_\alpha(d^\alpha_m\cap
\tau^\alpha_m): m<\omega\}$ and $\langle d^\alpha_m:
m<\omega\rangle$ is a maximal antichain.

If $\alpha\in J $ we choose
$\eta^\alpha,d^\alpha_n,\tau^\alpha_m$, satisfying the
above  and let $b^\alpha_m=h_\alpha(d^\alpha_m),
c^\alpha_m=h_\alpha(d^\alpha_m\cap \tau^\alpha_m)$.

So ``$\alpha\in J$'' means that I played Boolean Algebras
and cadomorphisms as is the previous remark and we get in
the Boolean Algebra with some properties.

The desired Boolean algebra $B$ is $B_\alpha^*$.  We
shall investigate it and eventually prove it is
endo-rigid (in 3.11) and indecomposable   (in 3.12)
(3.1(1)), 3.1(2) are trivial).

Note also

\proclaim{3.3 Fact}
\roster\item For $\nu \in {}^{\omega>}\lambda,x_\nu$ is
free over $\{x_n:\, \eta \in{}^{\omega>}\lambda,\eta
\not=\nu\}$ hence also over the subalgebra of $B^c_0$ of
those elements based on
$\{x_\eta:\,\eta\in{}^{\omega>}\lambda,\eta\not= \nu\}$.

\item  For every branch $\eta$ of $F^\alpha$  such that
$\eta\not=\eta_\beta$ for $\beta<\alpha,\xi(\alpha)$; and
finite $\omega\subseteq \alpha$ {\it there is} $k$
such that $\{\rho:\eta\restriction k\leq \rho \in T\}$
is disjoint to ${}^{W>}\xi \cup \cup\{N^\beta \cap
T:\beta \in w,\beta +2^{\aleph_0} \leq\alpha\} \cup
\cup \{\overline d(\tau^\beta_n):n<\omega\in w\}$
\endroster
From 3.2 we can conclude:\endproclaim
\bigskip
\proclaim{3.4 Fact}  If
$\xi<\zeta(\beta),\beta<\alpha,I\subseteq T$ finite then
every element of $B_\alpha$, based on
$I\cup{}^{\omega>}\xi$ is in $B_\beta$.
\endproclaim
\bigskip
\proclaim{3.5 Notation}
\roster
\item  Let  $B^\xi$ be the set of $a\in B^c_0$ supported
by ${}^{w>}\xi$
\item For $x\in B^c_0,\xi<\lambda$ let
$pr_\xi(x)=\cap\{a\in B^\xi:\,x\leq a\}$.
\item For $\xi <\lambda$ let
$\epsilon(\xi)=Min\{\gamma:\zeta(\gamma)>\xi\}$.

\item For $\gamma,\alpha^*$ let $B_{<\gamma>}=\langle
\{x_\eta:\,\eta\in{}^{\omega>}\zeta(\gamma)\}\cup
\{a_\beta:\beta<\gamma\}\rangle$.
\item For $\xi<\lambda$ let $B_{[\xi]}=\langle\{x_\eta:\,
\eta\in {}^{W>}\xi\}\cup\{a_\beta:\zeta(\beta)\geq
\xi\}\rangle_{B^c_0}$.
\endroster
\endproclaim
\bigskip
\proclaim{3.6 Fact}
\roster \item $B^\xi$ is a complete Boolean subalgebra of
$B^c_0$.
\item $pr_\xi(x)$ si well defined for $x \in B^c_0$
\item if $\xi_0<\xi_1<\lambda,x\in B^c_0$ then
$pr_{\xi_0}(pr_{\xi_1}(x)=pr_{\xi_0}(x)$.
\item If $\xi<\lambda,w\subseteq T$ is finite then
the function
$pr_{\xi,w}(x)=\cap\{y\in\langle B^\xi\cup\{x_\nu:\,
\nu\in w\}\rangle: x\geq y\} $ is well defined.
\endroster
\endproclaim
\bigskip\proclaim{3.7 Fact} \roster
\item For $x\in B_{\alpha^*},\xi<\lambda$, the element
$pr_\xi(x)$ belongs to $B_{[\xi]}$
\item for $x\in B_{\alpha^*},\xi,\lambda,w\subseteq
{}^{\omega>}(\xi+1)$, the element $pr_{\xi w}(x)$
belongs to $B({}^{\omega>}\xi,w)$.
\endroster
\endproclaim
\bigskip
\demo{Proof} 
 We prove this for $x\in B_\alpha$, by induction on
$\alpha$ (for all $\xi$)

Note that $pr_\xi( \cup_{\ell<n} x_\ell)=\cup_{\ell<n}
pr_\xi(x_\ell)$.
\bigskip
{\bf Case i}: $\alpha=0$, or even
$(\forall\beta,\alpha)[\zeta(\beta)\leq\xi]$

Easy; if
$x=\tau(_0,\dots,a_{n-1},x_{\nu_0},\dots,x_{nu_{m-1}})$
where $\tau$ is a Boolean term, $a-\ell\in
B_{[xi]},\nu_\ell \in {}^{\omega>}\xi$; by the remarks
above w.l.o.g $x=\cap_{\ell<n+m}\tau_\ell,
\tau_\ell\{a_\ell,1-a_\ell\}$ when
$\ell<n,\tau_\ell\in\{x_{\nu_{\ell-n}},1-x_{\nu_{\ell-n}}
\}$
when $n\leq n+m$, and the sequence $\langle x_{\nu_0}
,\dots,x_{\nu_{n-1}}\rangle$ is with no repetition, then
clearly  $pr_\xi(x)=\cap_{\ell<n}\tau_\ell\in B_{[\xi]}$;
\bigskip
{\bf Case ii}: $\alpha$ limit.

Trivial as $B_\alpha=\cup_{\beta<\alpha}B_\beta$.
\bigskip
{\bf Case iii}: $\alpha=\beta+1$

By the induction hypothesis w.l.o.g. $x\not\in B_\beta$. 
As $x\in B_\alpha$ there are disjoint $e_0,e_1,e_2\in
B_\beta$ such that $x=e_0\cup(e_1\cap a_\beta) \cup
(e_2-a_\beta)$.  It suffices to prove that
$pr_\xi(e_0),pr_\xi(e_1\cap a_\beta),
pr_\xi(e_2-a_\beta)\in B_{[\xi]}$, the first is trivial
and w.l.o.g. we concentrate on the second.  There are
$\xi_0<\zeta(\beta)$ and $k<\omega$ such that $e_1$ is
based on $J\buildrel  def\over =^{\omega>}\lambda
-\{\rho:\eta_\beta\restriction k\leq \rho
\in{}^{\omega>}\lambda\}$ and each $d^\beta_n(n<\omega)$ is
based on ${}^{\omega>}\xi_0$.  By case i, we can assume
$\xi<\zeta(\beta)$ hence w.l.o.g. $\xi<\xi_0$, and by the
induction hypothesis and 3.6(3) it suffices to prove
$pr_{\xi_0}(e_1\cap a_\beta)\in B_{[\xi]}$.  W.l.o.g.
$e_1\cap d^\alpha_m=0 $ for $m<k$ and now clearly
$pr_{\xi_0}(e_1\cap a_\beta)=e_1$ as $pr_{\xi_0} (e_1\cap
d^\alpha_m\cap\tau^\alpha_m)=e_1\cap d^\alpha_m$ for
$m\geq k$, (because $d^\alpha_m, e_1$ are based on
$J,{}^{\omega>}\xi_0\subseteq J$ and $\tau^\alpha_m$ is
based on ${}^{\omega>}\lambda-J$ and is $>0)$

\indent 2) Same proof.
\bigskip
\proclaim{3.8 Lemma}  Suppose $I, w$ satisfy:

$(*)_{I,w} I\subseteq {}^{\omega>}\lambda,w \subseteq
a^*, I$ is closed under initial segments, and for every
$\alpha<\alpha^*$ {\it if}
$\wedge_{m<\omega}(\eta_\alpha\restriction m\in I)$ then
$\tau^\alpha_m, d^\alpha_m$ are based on $i$ and belong
to $B(I<w)$.

Then for any countable $C\subseteq B_\alpha^*$ there is a
projection from $\langle B(i,w),C\rangle_{B_0^c}$
onto $B(I,w)$.
\endproclaim
\bigskip
\demo{Proof}  We can easily find $I(*),w(*)$ such
that $C\subseteq W(*)$, $w\subseteq
w(*)\subseteq \alpha^*, |w (*)-w|\leq
\aleph_0, I\subseteq I(*)\subseteq {}^{\omega>}\lambda,
|I(*)-I|\leq \aleph_0$ and if $\alpha \in
w(*)-w$, then $\tau^\alpha_m,d^\alpha_m\in
B(I*),w(*))$.  Let
$w(*)-w=\{\alpha_\ell:\ell<\omega\}$, and we define
by induction on $\ell$ a natural  number $k_\ell<\omega$,
such that the sets $\{\nu\in{}^{\omega>}\lambda:\nu$
appears is $\tau^{\alpha_\ell}_m$ for some $m>k_\ell\}$ are
pairwise disjoint to $I$.  Now we can extend the identity
on $B(I,w)$ to a projection $h_0$ from
$B(I(*),w)$ onto $B(I,w)$ such that of
$\ell<\omega,m>k_\ell$, then
$h_0(\tau^{\alpha_\ell}_m\cap d^{\alpha_\ell}_m)=0$.  Now
we can define by induction on $\alpha\in
(w(*)-w)\cup \{0,\lambda\}$ a projection
$h_\alpha$ from $B(I(*),w \cup (w(*) \cap
\alpha))$ onto $B(I,w)$ extending $h_\beta$ for
$\beta<\alpha$ and $\beta\in (w(*)-w\cup
\{0\})$.  For $\alpha=0$ we have it defined, for
$\alpha=\lambda$ we get the conclusion, and in limit
stages  take the union.  In successive stages there is no
problem by the choice of $h_0$, and the
$k_\ell$'s).
\enddemo
\bigskip
\proclaim{3.9 Claim}   If $B'$ is an uncountable
subalgebra  of $B_{\alpha'}$ {\it then} there is an
antichain $\{d_n:n<\omega\}\subseteq B'$ and for no $x\in
B, x \cap d_{2n}=0,x\cap d_{2n+1}=d_{n+1}$ for every $n$
provided that

(*) no one countable $I\subseteq {}^{\omega >}\lambda$ is
a support for every $a\in B'$.
\endproclaim
\bigskip
\demo{Proof}  We now define by induction on
$\alpha<\omega_1,d_\alpha,I_\alpha$, such that:

(i)  $I_\alpha\subseteq {}^{\omega>} \lambda$ is
countable.

(ii) $\cup_{\beta<\alpha} I_\beta\subseteq I_\alpha$ and
for  $\alpha$ limit, equality holds.

(iii) $d_\alpha \in B'$ is supported by $I_{\alpha+1}$

There is no problem with this.

By (iii) for each $\alpha$ there are $ \tau^0_\alpha\in
\langle a_\eta:\eta\in I_\alpha\rangle_{B^c_0},
\tau^1_\alpha, \tau^2_\alpha \in \langle a_\eta:\eta\in
I_{\alpha+1}-I_\alpha\rangle_{B^c_0}$ such that
$\tau^1_\alpha \cap \tau^2_\alpha=0,\,
\tau^0_\alpha\cap\tau^1_\alpha \leq d_\alpha$,
$\tau^0_\alpha \cap \tau^2_\alpha \leq 1-d_\alpha$.

By Fodor's lemma w.o.l.g. $\tau^0_\alpha=\tau^0$ (i.e. does
not depend on $\alpha$).  For each $\alpha$ there is
$n(\alpha)<\omega$ such that

$\tau^0_\alpha\in \langle a_\eta:\eta\in I_\alpha \cap
{}^{n(\alpha)\geq}\lambda \rangle_{B^
c_0}, \tau^1_\alpha,
\tau^2_\alpha \in \langle a_\eta:\eta\in (I_{\alpha
+1}-I_\alpha)\cap
{}^{n(\alpha)\geq}\lambda\rangle_{B^c_0}$ 

Again by renaming w.o.l.g $n(\alpha)=n(*)$ for every
$\alpha$.  Let for $n<\omega, d^n=d_n-\cup_{\ell<n}
d_\ell,\tau^n=\tau^0\cap\cap_{\ell<n}
\tau^2_\ell\cap\tau^1_n$, so easily $D^n\in
B',\langle d^n:n<\omega\rangle$ is an antichain,
$\tau^n\leq d^n$ and $\tau^n\langle a:\eta:\eta\in
{}^{n(*)\geq}\lambda\rangle_{B^c_0}$.  Suppose $x\in B',
x\cap d^{2n}=0, x\cap d^{2n+1}=d^{2n+1}$.  Hence for
$n<\omega, x\cap \tau^{2n}=0, x\cap
\tau^{2n+1}=\tau^{2n+1}$.  but by 3.8 (for
$I={}^{n(*)\geq}\lambda \omega=\phi$ and $c=\{x\})$,
there is such $x$ in $\langle a_\eta;\eta\in{}^{n(*)\geq}
\lambda\rangle_{B^c_o}$, an easy contradiction.

So we have proven that for every $\aleph_1$-compact
$B'\subseteq B_{\alpha^*}$, some countable $I\subseteq
{}^{\omega >}\lambda$ support every $x\in B'$.
\enddemo
\bigskip
\proclaim{3.10 Claim}  No infinite subalgebra $B'$ of
$B_{\alpha^*}$ is $\aleph_1$-compact.
\endproclaim
\bigskip
\demo{Proof}  Suppose there is such $B'$, and let $\xi$
be minimal such that there is such $B'\subseteq B_{[\xi]}$
\enddemo
\bigskip
{\bf Part I}  if (*) 
\roster
\item"(a)" $B'\subseteq B_{\alpha^*}$ is
$\aleph_1$-compact and infinite and

\item"(b)" $B'\subseteq B_{[\xi]}$

{\it then}

\item"(c)"  for every $\zeta<\xi$ and $x\in B'-\{y:\{z\in
B':z\leq y\}$ is fintie $\}$, there is $x_1\in B',
x-1\leq x$ such that for no $y\in B_{[\zeta]}, y \cap
x=x_1$.

So assume $B'$, satisfies a) and b) but they fail c) for
$\zeta<\xi$ and $x\in B'$, where $\{y:y\leq x,y\in B'\}$
is finite.  So for every $z\in B'$, there is $g(z)\in
B_{[\zeta]}$ such that $g(z)\cap x=z\cap x$ (use
$x_1=z\cap x$).  Let $B^a$ be the subalgebra of
$B_{[\zeta]}$ generated by $\{g(z):z \in B'\}$.  Clearly 
$\{y\in B':y \leq x\}=\{t\cap x:t\in B^a\}$.  Let
$x^*=pr_\zeta(x)$, (it is in $B_{[\zeta]}$ by 3.7(1)) and
let $ B^b=\{t\cap x^*;t \in B^a\}\cup \{t \cup
(1-x^*):t\in B^a\}$.  Clearly $B^b$ is a subalgebra of
$B_{[\zeta]}$, and $1-x^*$ is an atom of $B^b; B^b$ is
infinite as there are in $B'$ distinct $x_n\leq x$, so
$g(x_n)\in B^a$ hence $g(x_n)\cap x^*\in B^b$.  As $x\leq
x^*$ and $[n\not=\Rightarrow g(x_n) \cap x\not= g(x_m)\cap
x]$ clearly
 $[n\not= m\Rightarrow g(x_n)\cap x^* \not- g(x_n)\cap
x^*]$.  We shall prove that $B^b$ is $\aleph_1$-compact,
thus contradicting the choice of $\zeta$.  Let $d_n
\in B^b$ be pairwise  disjoint, and we want to find $T\in
B^b, t\cap d_{2n}=0, t \cap d_{2n+1}= d_{2n+1}$ (for
$n<\omega)$.  Clearly w.o.l.g. $d_n \leq x^*$ (as $1-x^*$
is an atom of $B^b$).  So $d_n=t_n\cap x^*$ for some
$t_n\in B^a$, hence easily $t_n \cap x\in B'$ so for some
$x_n\in B', x_n\leq x$ and $t_n\cap x=x_n\cap x=x_n$.  So
$X-n=g(x_n)\cap x$.

For $n\not=m, x_n\cap x_m=(t_n\cap x)\cap (t_m\cap x)\leq
(t_n\cap x^*)\cap (t_m\cap x^*)= d)n\cap d_m=0$.

As $B'$ is $\aleph_1$-compact there is $y\in B', y \cap
x_{2n}=0, y\cap x_{2n+1}=x_{2n+1}$.  Now $g(y), d_n, t_N$
belongs to $B_{[\xi]}$ and (as $x_n\leq x \leq x^*)$:
\smallskip
\ \ \ \ (i)  $g(y)\cap d_2n\cap x=g(y)\cap t_{2n} \cap x=
g(y) \cap x_{2n} \cap x=y\cap x_{2n} \cap x=0$.

\ \ \ \ (ii)  $(y)\cap d_{2n+1} \cap x=g(y) \cap
t_{2n+1}\cap x=g(y) \cap x_{2n+1} \cap x= y \cap
x_{2n+1}\cap x =x_{2n+1} \cap x=t_{2n+1}\cap x=d_{2n+1}
\cap x$.

Now by the definition of $X^*=pr_\zeta(x),[\tau \in
B_{[\zeta]}\wedge \tau \cap x=0\Rightarrow \tau\cap
x^*=0]$.  As $(1-\tau\in B_{[\zeta]}, x\leq 1-\tau)$ hence
by (i) (for $\tau=g(y)\cap d_{2n}$):

\ \ \ \ (iii) $g(y)\cap d_{2n}\cap x^*=0$

Also by the definition of $x^*=pr_\zeta(x)$:

$\tau_1,\tau_2\in B_{[\zeta]}\wedge \tau_1\cap x=
\tau_2\cap x\Rightarrow x^*=\tau_2\cap x^*$ (as
$\tau_1-\tau_2\in B_{[\zeta]}, x\leq 1-(\tau_1-\tau_2))$
hence by (ii)

\ \ \ \ (iv)  $g(y)\cap d_{2n+1} \cap x^*=d_{2n+1} \cap
x^*$.

But $d_n\leq x*$, so from (iii) and (iv) $(g(y) \cap
d_{2n}=0, (g(y)\cap x^*)\cap d_{2n+1}=d_{2n+1}$, and
$g(y)\in B^a$ hence $g(y) \cap x^*\in B^b$.  So $B^b$ is
$\aleph_1$-compact this contradicits the minimality of
$\xi$, so we finish Part I.
\endroster
\bigskip
\bigskip
{\bf Part II}:  If $B^1$ is $\aleph_1$-compact
$B^1\subseteq B^2, B^2=\langle B^1\cup\{z\}\rangle$ then
$B^2$ is $\aleph_1$-compact.
\bigskip
The proof is straightforward.  [If $d_n\in B^2$ are
pairwise disjoint, let $d_n=(d^1_n\cap z)\cup (d^2_n-z)$
for some $d^1_n, d^2_n\in B^1$.  Now w.l.o.g. $D^1_n\cap
d^1_m=0$ for $n\not=m$- otherwise replace then by
$d^1_n-\cup_{\ell<n} d^1_\ell$; similarly $d^2_n\cap
d^2_m=0$, for $\not=m$.  So there are $y^\ell\in B^1,
y^\ell\cap
d^\ell_{2n+1}=d^\ell_{2n+1}=d^\ell_{2n+1}$, and $(y^1\cap
z) \cup (y^2-z)$ is the solution.]
\bigskip
\bigskip
{\bf Part III}:  $\xi$ cannot be a successor ordinal.
\bigskip
\demo{Proof}:  Let $B'$ satisfy (*)

Suppose $\xi=\zeta+1$, and by 3.9 there is a countable
$I\subseteq {}^{\omega>}\xi$ which support every $a\in
B'$.  w.l.o.g.$I$ is closed under initial segments and
$k=|I-{}^{\omega>\zeta}$ is minimal.  Now part I can be
applied with $\langle B_{[\zeta]},\{a_\eta:\eta\in w\}
\rangle_{B^c_0}$, for any finite $w\subseteq I$ of power
$< k$ instead $B_{[\zeta]}$ (using 3.7(2) instead
3.7(1)).  So by applying Part I (to $\langle
B_{[\zeta]},\{a_\eta:\eta\in w\}\rangle_{B^c_0}$) we can
add to its conclusion:

d)  {\it for every} finite $w\subseteq
I,|w|<|I-{}^{\omega>}\zeta|$ and $x\in B'$ and $x\in B'$
for which $\{y\in B':y\leq x\}$ is infinite, {\it there
is} $x_1\in B',x_1\leq x$ {\it such that for} no  $y\in
\langle B_{[\zeta]}\cup \{a_\eta:\eta\in
w\}\rangle_{B^c_0}, y \cap x=x_1$.

Now $I-{}^{\omega>}\zeta$ is infinite [otherwise let
$B''= \cup\{a_\eta;\eta \in
I-{}^{\omega>}\zeta\}\rangle_{B^c_o}$, easily it is
infinite and $\aleph_1$-compact by Part II and then we
apply Part I: for

$I_{}^{\omega>}\zeta=\{\eta_0,\dots,\eta_{k-1}\}$ and
for $u\subseteq\{0,\dots,k-1\}$, let
$x_u\buildrel def\over= \cap \{x_{\eta_\ell}:\ell\in
u\} \cap \{1-x_{x\ell}:\ell<k, \ell\not\in u\}$ so $x_u\in
B",1=\cup\{x_u:u\subseteq \{0,\dots,k-1\}\}$, hence for
some $u,\{y\in B':y\leq x_u\}$ is infinite; $\zeta,x_u$
contradict the conclusion of Part I.

As $B'$ is $\aleph_1$-compact, for any $x\in B'$ such
that $\{y\in \leq x\}$ is infinite, $x$ can be splitted
in $B'$ to two elements satisfying the same i.e.
$x=x^1\cup x^2, x^1\cap x^2=0,\{y\in B':y\leq x^\ell\}$
is infinite for the $\ell=1,2$.  Let
$I_{}^{\omega>}\zeta=\{\eta_\ell:\ell<\omega\}$, so we
can find pairwise disjoint $e_n\in B'$, such that
$e_n=d_{2n} \cup d_{2n+1}, d_{2n} \cap d_{2n+1}=0$ and
that for no $y\in \langle B_{[B_\zeta]}\cup
\{a_{n\ell}:\ell<n\}\rangle, y\cap (d_{2n} \cup d_{2n+1})=
d_{2n+1}$ for every $n$.  So for no $n\, y\in \langle
B{[\zeta]}\cup \{a_{n\ell}:\ell<n\}\rangle_{B^c_0}$.

As $y\in B'$ clearly $y\in B_{[\zeta+1]}$, but $y$ is
based on ${\omega>}\zeta\cup\{a_{\eta\ell}:\ell<\omega\}$
so by 3.7(2) $y \in \langle B_{[\zeta]}\cup
\{a_{\eta\ell}:\ell<\omega\}\rangle_{B^c_0}$,
contradiction to $y \cap(d_{2n}\cup d_{2n+1})0=d_{2n+1}$.
\enddemo
\bigskip
\bigskip
{\bf Part IV}:  Let $B'$, satisfy (*) of Part I.  by 3.9
for some countable $I\subseteq {}^{\omega>}\xi$, every
$b\in B'$ is based on $I$.  By Part III $\xi$ is not a
succesor ordinal, so necessarily $cf(\xi)=\aleph_0$, let
$Fi(B')=\{x\in B':\{y\in B':y\leq x\}$ is finite$\}$. 
Next we shall show:

(**)  for some finite
$w\subseteq\{\gamma:\zeta(\gamma)=\xi\}$ and $x^* \in
B'-Fi(B')$ for every $y<x^*$ from $B'$, for some $Z\in
\langle \cup_{\zeta<\xi}
B{[\zeta]}\cup\{a_\alpha:\alpha\in w\}\rangle_{B^c_0},
z\cap x^*=y$.

Suppose (**) fails and we define by induction on
$n<\omega, x_n,y_n,w_n$, such that:
\roster
\item"(i)" $x_n\in B'$,
\item"(ii)" $1-\cup_{i<n} x_i \not\in Fi(B')$
\item"(iii)" $w_n\subseteq \{\gamma:\zeta(\gamma)=\xi\}$
is finite
\item"(iv)" $w_n\subseteq w_{n+1}$
\item"(v)" $y_n\leq x_n,y_n\in B'$.
\item"(vi)" for no $z\in \langle \cup_{\zeta<\xi}
B_{[\zeta]}\cup\{a_\alpha\in w_n\}_{B^c_0}$ is $\cap
x_n=y_n$.
\endroster
For $n=0\, 1 \not\in Fi(B')$

For every $n$ let $w_n$ be a finite subset of
$\{\gamma:\zeta(\gamma)=\xi\}$ extending $\cup_{\ell<n}
w_\ell$, such that for every $\ell<n, x_\ell,y_\ell\in
\langle\cup_{\zeta<\xi} B_{[\zeta]} \cup
\{a_\alpha:\alpha \in w_n\}\rangle_{B^c_0}$.  Then as
$1-\cup_{\ell<n} x_i \not\in Fi(B')$, and as $B'$ is
$\aleph_1$-compact, there is $x_n\leq 1-\cup_{i<n} x_i,
x_n\in B', 1-\cup_{\ell\leq n} x_i\not\in Fi(B')$ and
$x_n\not\in Fi(B')$.  Now as (**) fails $w_n,x_n$ does
not satisfy the requirements on $w,x^*$ in (**), so
there is $y_n\in B', y_n\leq x_n$ such that for no $z\in
\langle\cup_{\zeta<\xi} B_{[\zeta]}\cup
\{a_\alpha:\alpha \in w_n\}\rangle_{B^c_0}$.  By 3.8 we
can easily get a contradicition to (vi).  So (**) holds.

Let $t_0,\dots,t_m \in B_{[\xi]}$ be such that
$\cup^m_{\ell=1} t_\ell=1$ and $(\forall \ell \leq m)
(\forall \alpha \in w)\, [t_\ell\leq a_\alpha \vee t-\ell
\cap a_\alpha=0]$.  There is an $\ell\leq m$ such that
$\{u\cap t_\ell:y\leq x^*$ and $y\in B'\}$ is finite.  It
is clear (by Part II) that $B''=\langle
b',t_\ell\rangle_{B^c_0}$ is $\aleph_1$-compact:  also
$X^*\cap t_\ell \in B''-Fi)B'')$.  Now if $y\in B'', y\leq
x^*\cap t_\ell$ then for some $y'\in B',y=y'\cap t_\ell$
and w.l.o.g. $y'\leq x^*$, so for some $z\in \langle
\cup_{\zeta<\xi} B_{[\zeta]}\cup \{a_\alpha:\alpha \in
w\}\rangle_{B^c_0} z\cap x^*=y'$ hence $z\cap (x^*\cap
t_\ell)=y$, and by the choice of $t_\ell$, for some
$z'\in \cup_{\zeta<\xi} B_{[\zeta]}$, the equation
$z'\cap (x^*\cap t_\ell)=z\cap (x^*\cap t_\ell)=y$ holds.

So $B'',x**\buildrel\text{def}\over = x^* \cap t-\ell$
satisfy requirements in (**).  Now we use (c) of part I. 
As $cf(\xi)=\aleph_0$, let $\xi=\cup_{n<\omega}
\zeta_n$, and we define by induction on
$n<\omega,x_n,y_n$ such that:
\roster
\item"(i)" $x_n\in B'', x_n\leq x^{**}$
\item"(ii)" $x^{**} -\cup_{\ell<n} x_i\not\in Fi(B'')$
\item"(iii)" $y_n\in B', y_n\leq x_n$
\item"(iv)" for no $z \in B_{[\zeta_n]}, z\cap x_n=y_n$.
\endroster

As $B"$ is $\aleph_1$-compact, for some $x^*\in B",
z^*\cap x_n=y_n$ for each $n$.

Now as $B'',x^{**}$ satisfy (**), for some $z^{**}\in
\cup_{\zeta<\xi} B_{[\zeta]} z^*\cap x^{**}=z^{**}\cap
x^{**}$.  So for some $n\, z^{**}\in B_{[\zeta_n]}$,
contradicting (iv) above.  Thus we have finished the proof
of 3.10. \bigskip
\proclaim{3.11 Claim}
$B_{\alpha^*}$ is endo-rigid.
\endproclaim
\bigskip
\demo{Proof} Suppose $h$ is a counterexample, i.e. $h$ is
an endomorphism of $B_{\alpha^*}$ but $B_{\alpha^*}/\,
Ex\, Ker(h)$ is infinite, and we shall get a
contradiction.

Clearly if for some $\alpha, N^\alpha=(B^*,h)\, h$
maps $B^*\cap B_{\alpha^*}$, into itself and $\alpha\in J$
(see Stage B) then $h(a_\alpha)$ realizes the type
$p_\alpha$, contradiction (by stage A, $B_{\alpha^*}$
omits $p_\alpha$).  So we shall try to find such $\alpha$
which satisfies the requirements in Stage B for belonging
to $J$.  We assume $N^\alpha=(|N^\alpha|,h_\alpha),
|N_\alpha|\subseteq B_\alpha, h_\alpha=h \restriction
N^\alpha, h_\alpha$ maps $N^\alpha\cap B_\alpha$ onto
itself, and $N^\alpha_0$ contains some elements we need
and somewhat more (see latter).  As $W$ is a barrier this
is posible.  We then will choose $\eta_\alpha$, an
$\omega$-branch of $f^\alpha$, distinct from $\eta_\beta$
for $\beta<\alpha$ [if $\beta+2^{\aleph_0} \leq \alpha$
this follows, the rest exclude $<2^{\aleph_0}$ branches
of $f^\alpha$ but there are $2^{\aleph_0}$ such
branches], a maximal antichain $\langle
d_n: n<\omega\rangle$ of $B_\alpha,d_n\in N^\alpha_0$, and
$\tau_n\in N^\alpha$ in $\langle
x_\nu:\eta\alpha\restriction n<\nu\in T\rangle_{B^c_)}$,
and let
$b_n=h(d_n),c_n=h(d_n\cap\tau_n),p_\alpha=\{x\cap
b_n=c_n:n<\omega\}$, and $a_\alpha=\cup_{n<\omega}
(d_n\cap \tau_n)\in B^c_0$.  All should have superscript
$\overline d,\overline\tau$ (where $\overline
d=\langle:n<\omega\rangle,\overline \tau=\langle
\tau_n:n<\omega\rangle)$ but we usually omit  them or
write $a_\alpha[\overline\tau,\overline d],
p_\alpha[\overline\tau\overline d]$ etc.

The choice of $\overline d,\overline \tau$ ( and
$\eta_\alpha$ which is determined by $\overline\tau$) is
done by listing the demands on them (see Stage B) and
showing a solution a solution exists.  The only
problematic one is (a) (omitting $p_\beta$ for $\beta\leq
\alpha$) and we partition it to three cases:
\roster
\item"(I)" $\zeta(\beta)<\zeta(\alpha)$ or
$\zeta(\beta)=\zeta(\alpha),\beta+2^{\aleph_0} \leq
\alpha$,
\item"(II)"
$\zeta(\beta)=\zeta(\alpha),\beta<\alpha<\beta+2^{\aleph_0}$.
\item"(III)" $\beta=\alpha$
\endroster

We shall prove that every $\overline \tau,\overline d$ are
O.K. for (I), that for any family $\{(\overline d^i,
\eta^i,\overline\tau^i): i<2^{\aleph_0}\}$ ($\eta$ a branch
of $f^\alpha$, etc.) with pairwise distinct $\eta^i$'s, all
except $2^{\aleph_0}$ many are O.K. for instance of (II),
and that there is a family of $2^{\aleph_0}$ triples
$(\overline d,\overline \eta,\overline\tau)$ satisfying
(III) with pairwise distinct $\eta^i$'s.  This clearly
suffices.
\bigskip
{\bf Case I}: $\zeta(\beta)<\zeta(\alpha)$ or
$(\zeta(\beta)=\zeta(\alpha),\beta+2^{\aleph_0} \leq
\alpha$

Suppose some $x\in\langle_\alpha,a_\alpha[\overline
\tau,\overline d]\rangle_{B^c_0}$ realizes $p_\beta$. 
Clearly there is a partition $\langle e_\ell:\ell <
4\rangle$ of 1 (in $\beta_\alpha$) such that
$x=e_0\cup(e_1\cap a_\alpha [\overline \tau,\overline
d])\cup (e_2-a_\alpha [\overline \tau,\overline d])$. 
Choose $\xi<\zeta(\alpha)$ large enough and finite
$w\subseteq \alpha$ so that
$[\zeta(\beta)<\zeta(\alpha)\Rightarrow
\zeta(\beta)<\xi], d_n,h_\alpha(d_n) b^\beta_n$, are
based on $\{x_\nu:\nu\in {}^{\omega>}\xi\}$ (for
$\eta<\omega)$ and
$c^\beta_\ell(\ell<\omega),e_0,e_1,e_2,e_3$ are based on
$J=\{\nu\in T:\eta_\alpha\restriction k\not\leq \nu\}$,
where $k<\omega$ also satisfies such that
$\eta_\alpha(k)>\xi,\eta_\alpha k\not\in N_\beta$.

We claim:

(*) there is $m<\omega$ such that $b^\beta_m\cap
(e_1\cup e_2)-\cup_{n\leq k} d_n\not= 0$

For suppose (*) fail, then as $a_\alpha[\overline
\tau,\overline d]\cap(\cup_{n\leq k} d_n)\in B_\alpha$,
w.l.o.g. $(e-1\cup e_2) \cap \cup_{n\leq k} d_n=0$
(otherwise let

\indent $e'_0=e_0\cup(e_1\cap a_\alpha[\overline
\tau,\overline d]\cap \cup_{n\leq k} d_n)\cup
(e-2\cap\cup_{n\leq k}
d_n-a_\alpha[\overline\tau,\overline d])$

\indent $e'_1=e_1-\cup_{n\leq k} d_n$,

\indent $e'_2=e_2-\cup_{n\leq k} d_n)$.

So for every $M<W$, $b^\beta_M\cap (e_1\cup e_2)=0$.

So if $x$ realizes $p_\beta$ then so does $e_0$, but
$e_0\in B_\alpha$ contradicting an induction hypothesis. 
So (*) holds.

Now as $\langle d_n:n<\omega\rangle$ is a maximal 
antichain in $B_\alpha$, for some $\ell<\omega ,
d_\ell\cap(b^\beta_m\cap(e_1\cup e_2-\cup_{n\leq k}
d_n))\not=0$.  Necessarily $\ell> k$.  So for some
$\epsilon \in\{1,2\},d_\ell\cap b^\beta_m\cap
e_\epsilon\not= 0$.  As $x$ realizes
$p_\beta,x\cap(d_\ell\cap b^\beta_m\cap
e_\epsilon)=d-\ell\cap c^\beta_n\cap e_\epsilon$ which is
based on $J$.  But we know that $x\cap (d_\ell\cap
b^\beta_m \cap e_\epsilon)$ is $d_\ell\cap b^\beta_m\cap
e_1\cap a_\alpha [\overline \tau,\overline d]=d_\ell \cap
b^\beta_m\cap e_1\cap\tau_\ell$ (if $\epsilon=1$) or
$d_\ell\cap b^\beta_m\cap e_2\cap (1-a_\alpha([\overline
\tau,\overline d])=d_\ell\cap b^\beta_m\cap e_2\cap
1-\tau_\ell$ (if $\epsilon=2$).  As $d_\ell\cap
b^\beta_m\cap e-z\not=0$ is based on
$J,\ell>k,\eta_\alpha(k)>\xi,\tau_\ell$ is free over $J$,
(see  Fact(2)).  Necessarily $x\cap (d_\ell\cap
b^\beta_m\cap e_z$) is not based on $J$, contradiction.
\bigskip
{\bf Case II}: $\beta<\alpha<\beta+2^{\aleph_0}$

We shall prove that if $\eta^\ell,\overline \tau^\ell$
are appropiate (for $\ell=1,2)$ and $\eta^1\not= \eta^2$
then $p_\beta$ cannot be realized in both
$\langle B_\alpha,a[\overline\tau_\ell,\overline
d]\rangle_{B^c_0}$.  (So as
$\beta<\alpha<\beta+2^{\aleph_0}$, there are less then
$2^{\aleph_0}$ non appropiate $\overline
\eta^?,\overline\tau^1)$.

As there is a perfect set of appropiate $\eta$'s it will
suffice to prove that for each $\omega$-branch $\eta$ of
Rang$(f^\alpha)$ for some appropiate $\tau_1\langle
B_\alpha,a^{\overline\tau}\rangle_{B^c_0}$ omits
$p_\alpha=p_\alpha[\overline\tau,\overline d]$ which will
be done in Case III.

Note that $I^\alpha_\beta=\{e\in B_\alpha:$ for some
$x\leq e$ for every $n\, x\cap b^n_\beta\cap e=
c^n_\beta\cap e\}$ is an ideal.

The details are easy.
\bigskip
{\bf Case III}: $\beta=\alpha$

This case is splitted into several subcases.  Let
$\eta_\alpha$ be any $\omega$-branch of
$f^\alpha,\eta_\alpha\not= \eta_\beta$ whenever
$\beta<\alpha<\beta+2^{\aleph_0}$.  Let 
$I^*=\cup\{\underline d(h(x));x\in B_\alpha\}$.  We shall
assume that $|I^*|\leq \aleph_0\Rightarrow I^* \subseteq
N^\alpha_0$, so in this case $p_\alpha$ is omitted by
$B_{\alpha+1}$ or $B_{\alpha^*}$ iff  ?? omitted by
$B_\alpha$ (by 3.7(1)).  As acomplishing this aim is
easier we shall  ?? this case (work as in III 4 and use
quite arbitrary $p_\beta)$.
\bigskip
{\bf Subcase III 1.}:  For some $\rho^*\in T$, and $a^*\in
B_\alpha-Ex\, Ker^*(h)$ for every $\leq \rho\in T$ for some
$\tau \in\langle x_\eta:\rho<\eta\in T\rangle_{B^
c_0},\tau a \cap a*\not=0=h(\tau\cap a^*)$. 

As we are interested not in $f^\alpha,N^\alpha)$ itself,
but in $h$, by using $Gm'(W)$, w.l.o.g. $\rho^*\in Range\,
(f^\alpha)$.  By 3.10 (for $rang\,(h)$, which by
assumption, is infinite) ?? easy manipulations (see 2.4
and [Sh 2]) there is maximal antichain
$??:n<\omega\rangle$ of $B_{\alpha^*}$ such that for no
$x\in B_\alpha$, $x\cap h(2_{2n}=h(d_{2n})$ and $\cap
h(d_{2n+1})=O$ W.l.o.g.$\{d_n:n<\omega\}\subseteq
N^\alpha_0$

It suffices to prove the conclusion for any
$\omega$-branch $\eta_\alpha$ of
$Range(f^\alpha),\rho^*<\eta_\alpha\not\in
\{\eta_\beta:\beta<\alpha\}$.  We define by induction on
$n,\tau_n\in N^\alpha_n< in
\langle x_\eta:\eta_\alpha\restriction n \leq
\eta\rangle_{B^c_0}, \tau_n \not=0,1$ and
$h(\tau_{2n})=1,h(\tau_{2n+1})=0$.  (possible by the
assumption of subcase III 1), so we finish this subcase.
\bigskip
{\bf Subcase III 2}.  For some $a^*\in
B_\alpha,\{h(x)-a^*:x\in B_\alpha, x\leq a^*\}$ is
infinite.

Clearly $ B^a=\{h(x)-a^*:  x\in B_{\alpha^*} x\leq
a^*\}\cup \{1-(h(x)-a^*): b_{\alpha^*}x\leq a^*\}$ is a
subalgebra of $B_{\alpha^*}$ (with $a^*$ an atom).  By
assumption (of ?? subcase) $B^a$ is infinite.  So by 3.9
there are $e_n\in B^a$, pairwise disjoint, and (? $x\in
B_a)\bigwedge_n (x\geq e_{2n} \wedge x\cap
e_{2n+1}=0)$.  As $a^*$ is an atom of $B^a$ w.l.o.g.
$\leq 1-a^*$, hence there is $d_n\leq a^*$  (in
$B_{\alpha^*}$, such that $h(d_n)=e_n$.  Clearly
$??-\cup_{\ell<n} d_n)=e_n-\cup_{\ell<n} e_\ell=e_n$, so
w.l.o.g. the $D_n$ are pairwise disjoint.  So by easy
manipulation for some $\langle d_n:n<w\rangle$ the
following holds: \roster
\item"(i)" $d_0=1-a^*$
\item"(ii)" $\langle d_n:n<\omega\rangle$ is a maximal
antichain of $B_{\alpha^*}$.
\item"(iii)"  for no $x\leq 1-a^* x\cap
h(d_{2n+2})-a^*$,  $x\cap h(d_{2n+1})-a^*=0$ \endroster

We can assume that $d_n,h(d_n)\in N^\alpha_0$.

Let $\overline \tau=\langle\tau^0_n:n<\omega\rangle$ be
a suitable suquence, (for our $\eta_\alpha$) then so are
$\overline \tau^\ell=\langle
\overline\tau^\ell_n,\omega\rangle$, for $\ell<4$ where:

\indent
$\tau^1_{2n}=1-\tau^0_{2n},\tau^1_{2n+1}=\tau^0_{2n+1}$;

\indent $\tau^2_{2n}=\tau^0_{2n},
\tau^2_{2n+1}=1-\tau^0_{2n+1}$;

\indent
$\tau^3_{2n}=1-\tau^0_{2n},\tau^3{2n+1}=1-\tau^0_{2n+1}$

Suppose for each $\ell<4$, in
$\langle_\alpha,a_\alpha[\overline
\tau^\ell,\overline d]\rangle_{B^c_0}$ there is an element
$y^\ell$ which satisfies $y^\ell \cap
h(d^\ell)-a^*=h(\tau^\ell_n\cap d_n)-a^*$ for $1\leq
n<\omega$.  W.l.o.g. $y^\ell\leq 1-a^*=d_0$ hence
$y^\ell\in B_\alpha$.  Now $(y^0\cup y^1\cap (y^2\cup
y^3)\in B_\alpha$ contradict (iii) above.

\bigskip
{\bf  Subcase III 3.}  For some $a^*\in B_\alpha^*-Ex\,
Ker^*(h)$, and $\rho^*\in T$, for every $\rho,\rho\leq
\rho\in T$ there is $\tau \in \langle x_\nu;\rho \leq \nu \in
T\rangle_{B^c_0}$ such that $h(\tau\cap a^*)\cap
a^*=\tau\cap a^*$

Clearly the function $h':b_\alpha^*\restriction
a^*\rightarrow B_\alpha^*\restriction a^*$ defined by
$h'(x)=h(x)\cap a^*$ is an endomorphism; W.l.o.g. the
assumption of subcase III 2 fails hence $\{(x)-a^*:x\leq
a^*\}$ is finite, hence the range of $h'$ is finite (as
$a^*\not\in Ex\, Ker^*(h)$, so by 2.4 there is $x\leq a^*$
such that $h(x)\cap a^*-x\not=0$; we know that $\underline
d(x)$ is countable, hence for some $\rho^{**},\rho^*
\leq\rho^{**}\in T, \{\nu:\rho\leq\nu\in T\}$ is disjoint to
$\underline(a^*) \cup\underline d(x))$.  Now by the
hypothesis of subcase III 3 we can easily find $\tau_n\in
\langle x_\nu:\rho^{**}\leq \nu \in t\rangle_{B_0^c}$, with
pairwise disjoint $\overline d(\tau_n)$ and $h(\tau_n\cap
a^*)\cap a^*=\tau_n\cap a^*$.  So

$h(\tau_n(\cap x)\cap (a^*-x)=h((\tau_n\cap a^*)\cap
x)\cap(a^*-x)=h(\tau_n\cap a^*)\cap h(x)\cap(a^*-x)=
(h(\tau_n\cap a^*)\cap a^*)\cap h(x)\cap
(a^*-x)=(\tau_n\cap a^*)\cap
h(x)\cap(a^*-x)=\tau_n\cap
h(x)\cap(a^*-x)=\tau_n\cap(h(x)\cap a^*-x)$

It is $\not=0$ [as $\overline d(\tau_n) \cap \overline
d(x) \cup \overline d(h(x)) \cup \overline
d(a^*))=\phi)$ and $h(x) \cap
a^*-x\not=0,\tau_n\not=0]$, and for different $n$ we get
different values.  So $\{h(y\cap x)\cap (a^*-x):x\in
B_\alpha'\}$, is finite.  Hence $\{(y\cap x)-x;y\in
B_{\alpha'}$ is infinite. Leading to the assumption of
subcase III 2(with $x$ here for $a^*$ there).
\bigskip
{\bf Subcase III. 4} For some $\rho^*\in T$, and $a^*\in
B_{\alpha^*}\,Ex\, Ker^*(h)$ for every $\tau\in
\langle_x\nu:\rho^* \not\leq \nu\in T\rangle_{B^c_0}
h(\tau\cap a^*)\cap a^*$ is based on $\{\nu:\rho^*\not\leq
\nu \in T\}$.

W.l.o.g. the hypothesis of subcase III 1 fails hence
$\{h(\tau\cap a^*):\tau\in \langle x_\nu:\rho^*\leq\nu\in
T\rangle_{B^c_0}\}$ is infinite.  As also w.l.o.g. the
hypothesis of subcase III 2 fails we get $\{(\tau\cap
a^*)\cap a^*:\tau \in \langle x_\nu:\rho^*\leq\nu \in
T\rangle_{B^c_0}\}$ is infinite.  So by 3.9 we can find
$d_n\in \langle x_\nu:\rho^*\leq \nu \in T\rangle_{B^c_0}$
such that $d_n:n<\omega\rangle$ is a maximal antichain in
$B^c_0$, and there is no $x\in B_{\alpha^*}, z\cap
h(d_{2n}=h(d_{2n})+h(d_{2n}), x\cap h(d_{2n+1})=0$, and
$d_0=1-a^*$.

As before we can assume $\rho^*\in Rang(f^\alpha)$ and
$d_n\in N^\alpha_0$ for $n<\omega$.  We suppose
$\eta_\alpha\not\in\{\eta_\beta: \beta <\alpha\}$ is an
$\omega$-branch of $f^\alpha,\rho^*\leq\eta_\alpha$.

For any suitable $\tau$ if $y[\overline \tau,\overline
d]\in \langle B_\alpha,a_\alpha[\overline \tau,\overline
d]\rangle_{B^c_0}$ satisfies $\tau_n\in \langle
x_\nu:\rho^*\leq \nu\in T\rangle_{B^c_0}$ and $y[\overline
\tau,\overline d]\cap h(d_n)$, (for every $n$) then by
3.3 we easily get $y[\overline\tau,\overline d]\in
B_\alpha$, and then get contradiction by trying four
$\tau$'s, as in subcase III 2.
\bigskip
{\bf Subcase III. 5.}  There are $\rho^*\in T$ and an
atomless countable subalgebra $Y\subseteq B_\alpha^*$ and
pairwise disjoint $c_\ell\in Y(\ell<\omega)$ such that
for every $\ell$ and $\rho_\ell\in\{\rho:\rho^*\leq \rho\in
T\}$ for some $\tau_\ell\in \langle x_\nu:\rho_\ell\leq\nu
\in T\rangle_{B^c_0}$, the following holds:  for no $x\in
B^c_0$ is $\overline d(x)\subseteq \{\nu:\rho_\ell \not\leq
\nu\in T\}$ and $z\cap h(c_\ell)\cap
c_\ell-\tau_\ell=h(c_\ell\cap\tau_\ell)\cap
c_\ell-\tau_\ell$.

Let $\langle D_n:n<\omega\rangle$ be a maximal antichain
of $B_{\alpha^*}$ such that $d_{2n}=c_{2n}$

So w.l.o.g. $Y\cup\{d_n:n<\omega\}\subseteq
N^\alpha_0,\rho^*\in Rang(f^\alpha)$ (using $Gm'(W))$, and
even $\rho^* <\eta_\alpha$, and each $N^\alpha_m$ is
closed under the functions $h$ and
$\rho_\ell\rightarrow_\ell$ (implicit in the assumption
of the subcase).

We can now choose by induction on $n,t_n\in N^\alpha_n$,

$$\tau_n\in\langle x-\nu:\eta_\alpha \restriction n\leq
\nu\in T\rangle_{B^c_0}$$

such that

(*)  (a) for even $n$, for no $x\in B^c_0$ based on
$\{\nu:\eta_\alpha\restriction n\not\leq \nu \in t\}$ is
$x\cap h(d_n)\cap d_n-\tau_n=h(d_n\cap\tau_n)\cap
d_n=\tau_n$.

Why is this sufficient?  We let $\overline
d=\langle d_n:n<\omega\rangle$ and $\overline
\tau=\langle\tau_n: n<\omega\rangle$.  So assume some
$y[\overline \tau,\overline d]\in \langle
B_\alpha,a_\alpha[\overline\tau,\overline
d]\rangle_{B^c_0}$ realizes $p_\alpha[\overline
\tau,\overline d]$, i.e. satisfies $y[\overline
\tau,\overline d]\cap h(d_n)=h(d_n\cap\tau_n)$ for every
$n$.  As $y[\overline \tau,\overline d]\in
\langle B_\alpha,a_\alpha[\overline \tau,\overline
d]\rangle_{B^c_0}$ for $y[\overline \tau,\overline
d]=e_0[\overline \tau,\overline d]\cup(e_1[\overline
\tau,\overline d]\cap a_\alpha[\overline\tau,\overline
d])\cup (e_2[\overline \tau,\overline
d]-a_\alpha[\overline\tau,\overline d])$.

For some $m(*),\omega, \underline d(e_0[\overline
\tau,\overline d])\cup \underline d(e_1[\overline
\tau,\overline d]\cup \underline d(e_2[\overline
\tau,\overline d])$ is disjoint to
$\{\nu:\eta_\alpha\restriction m(*) \leq \nu\in T\}$ ( see
3.3(2)).

Now we compute for $n$ even $>m(*)$:

$z\buildrel\text{def}\over = h(d_n\cap \tau_n)\cap
d_n-\tau_n=$

$=y[\overline\tau,\overline d]\cap h(d_n)\cap
d_n-\tau_n$ by the choice of $y[\overline\tau,\overline
d])$

$=(e_0[\overline \tau,\overline d]\cup (e_1[\overline
\tau\overline d]\cap a_\alpha[\overline \tau,\overline
d])\cup (e_2[\overline\tau,\overline d]-a[\overline
\tau,\overline d]))\cap h(d_n)\cap d_n-\tau_n=$

$=(e_0[\overline\tau,\overline d]\cap h(d_n)\cap
d_n-\tau_n)\cup ((e-1[\overline \tau,\overline d]\cap
a_\alpha[\overline\tau,\overline d])\cap h(d_n)\cap
d_n-\tau_n)\cup \cup((e_2[\overline\tau,\overline
d]-a_\alpha[\overline\tau,\overline d])\cap h(d_n)\cap
d_n-\tau_n)$

But $a_\alpha[\overline\tau,\overline d]\cap
d_n=\tau_n\cap d_n$ hence

$$(e_1[\overline\tau,\overline d]\cap
a_\alpha[\overline \tau,\overline d])\cap
d_n=(e_1[\overline \tau,\overline d\cap\tau_n\cap d_n$$

$$(e_2[\overline\tau,\overline d]-a_\alpha[\overline
\tau,\overline d])\cap d_n=(e_2[\overline \tau,\overline
d]-\tau_n)\cap d_n$$

Hence

$z=(e_0[\overline\tau,\overline d]\cap h(d_n)\cap
d_n-\tau_n)\cup (e_1[\overline\tau,\overline d]\cap
\tau_n) \cap h(d_n)\cap d_n-\tau_n)\cup
((e_2[\overline\tau,\overline d]-\tau_n\cap h(d_n)\cap
d_n-\tau_n)$

But the second term is zero and in the first $-\tau_n$ is
redundant, so 

$z=(e_0[\overline\tau,\overline d]\cap h(d_n)\cap
d_n-\tau_n)\cup e_2\cap h(d_n)\cap d_n-\tau_n)=$

$=(e_0[\overline \tau,\overline d]\cup e_2
[\overline\tau,\overline d])\cap h(d_n)\cap -\tau_n$

We can conclude

$(e_0[\overline \tau,\overline d] \cup e-2
[\overline\tau,\overline d])\cap h(d_n)\cap
d_n-\tau_n=h(d_n\cap\tau_n)\cap d_n-t_n$

contradicting the choice of $\tau_n$.

To finish Case III (hence the proof of 3(10) we need only
\bigskip
{\bf Why the five subcases exhaust all posibilities?}

Suppose none of III 1-5 occurrs.  By not subcase III 1
for some $\rho^0\in T$,
\roster
\item"(a)" $h(\tau)\not= 0$ for every $\tau\in \langle
x_\eta:\rho^0J\leq \eta\in T\rangle_{B^c_0}$

Let $Y$ be the $\langle
x_{\rho^0 \widehat{\, }<i>}:<\omega\rangle_{B_0^c}$.  As
$Y$ is countable, for some
$i(*)<\lambda,\{\nu:\rho^0\widehat{\,}<i(*)>\leq\nu \in
T\}$ is disjoint to $\cup\{\underline d(y)\cup\underline
d(h(y)):y\in Y\}$.  As `` not subcase III 5'' for some
$\rho^1,\rho^0 \widehat{\,}\langle i(*)\rangle \leq
\rho_1\in T$, and 
\item"(b)" there are no pairwise disjoint  non
zero $c_\ell\in Y(\ell<\omega)$, such that for every
$\rho^1_\ell,\rho^1 <\rho^1_\ell\in T$ for some
$\tau_\ell\in \langle x_\nu:\rho^1_\ell \leq\nu \in
T\rangle_{B^c_0}$, the following holds:

(*) for no $x\in B^c_0,\overline d(x) \subseteq
\{\nu:\rho^1_\ell \not\leq \nu \in T\}$
and $x\cap h(c_\ell)\cap c_\ell-\tau_\ell=h(c_\ell
\cap\tau_\ell)\cap c_\ell-\tau_\ell$

Clearly

\item"(c)" $\cup\{\underline d(y) \cup \underline
d(h(y)):Y\in Y\}$ is disjoint  to $\{\nu:\rho^1\leq\in T\}$

Let $Z+\{c\in Y$: for some $\rho^1_c,\rho^1\leq \rho^1_c\in
T$ for no $\tau\in
 \langle x_\nu:\rho^1
\leq\nu\in T\rangle_{B^c_0}$ does (*) of (b) hold  with
$c,\tau$ instead $c_\ell,\tau_\ell\}$

By (b) among any $\aleph_0$ pairwise disjoint members of
$Y$, al least one belong to $Z$.

It is quite easy to define $y_n\in Z ((n<\omega)$ such
that  $[y_n\in Ex\, Ker^*(h) \Rightarrow y_n\in Ex\,
Ker(h)]$, $[m<n\Rightarrow y_n\cap y_m=0]$, and  for every
$y\in Y-\{0\}$ for some $n,y\cap (\cup_{\ell<n}
y)\ell)\not= 0$ or $y_n\leq y$.  So (by the choice of $Y$)
$\langle y_n:n<\omega\rangle$ is maximal antichain of
${B^c_0}$.  We shall show $y_n\in Ex\, Ker(h)$; fix $n$
for a while, and suppose $y_n\not\in Ex\, Ker(h)$, and let
 $\rho^1_n,\rho^1_n\leq \rho^1_n\in T$ be such that for no
$\tau\in \langle x_\nu:\rho^1_n\leq\nu\in T\rangle_{B^c_0}$
does (*) of (b) hold.

Now for each $\tau\in \langle x_\nu:\rho^1_n\leq\nu\in
T\rangle_{B^c_0}$ as $y_n\in Z$, clearly [as (*) of (b)
fail for $y_n,\tau$ (and $\rho^1_n)]$ for some $x_1\in
B^c_0, \overline d(x_1)\subseteq \{\nu:\rho^1_n\not\leq \nu
\in T\}$ and $x_1\cap h(y_n)\cap
y_n-\tau=h(y_n\cap\tau)\cap y_n-\tau$.  Applying the
failure of (*) of (b) for $y_n,1-\tau,\rho^1_n$ we get
$x_2\in B^0_c,\overline d(x_2)\subseteq
\{\nu:rho^1_n\not\leq \nu\in t\}$ and $x_2\cap h(y_n)\cap
y_n-(1-\tau)=h(y_n\cap (1-\tau))\cap y_n(1-\tau)$; note
that $h(y_n\cap\tau)\leq h(y_n)$, and $h(y_n\cap
(1-\tau))=h(y_n)-h(y_n\cap\tau)$.  By these equations and
as $y_n h(y_n), x_1,x_2$ are based on
$\{\nu:\rho^1_n\not< \nu \in T\}$ (by (c) and their
choice resp.) clearly for some partition of $1,e_0^\tau,
e_1^\tau, e_2^\tau, e_3^\tau, \in B^c_0$, based on
$\{\nu:\rho^1\not\leq \nu\in T\}$:
 
\indent (i) $h(\tau\cap y_n)\cap y_n=e^\tau_0\cup
(e^\tau_1\cap\tau)\cup (e_2^\tau-\tau)$.

Now for any $\tau,\sigma\in \langle x_nu:\rho^1_n\leq\in
T\rangle$, easily (as $h$ is an endomorphism):

\indent (ii) $h((\tau\cup \sigma \cap y_n)\cap
y_n=(h(\tau\cap y_n)\cap y_n)\cap (h(\sigma \cap y_n)\cap
y_n)$.

\indent (iii) $h((\tau\cup \sigma \cap y_n)\cap
y_n=(h(\tau\cap y_n)\cap y_n)\cup (h(\sigma\cap y_n)\cap
y_n)$>

We can apply (i) to $\tau, \sigma$ and also to
$\tau\cup\sigma$, and substitute in (ii) (iii).

We get that 

\ \ \ \ ($\alpha$) $e^\tau_2\cap e^\sigma_2=0$ if
$\overline d(\tau)\cap \overline(\sigma)=),\tau,\sigma\in
\langle x_\nu;\rho^1\leq\in T\rangle_{b^c_0}$  (otherwise
substitute (i) in (ii) and intersect with $e_2^\tau\cap
e^\sigma_2)$ and get $(h((\tau\cap \sigma)\cap y_m)\cap
e^\tau_2\cap e^\sigma_2)= (e^\tau_2-\tau)\cap
(e^\sigma_2-\sigma)=e^\tau_2\cap \tau^\sigma_2\cap
(\tau\cup\sigma)$, and $(h((\tau\cap\sigma)\cap y_n)\cap
y_n) \cap (e^\tau_2\cap e^\sigma_2)\not\in \langle
\{x:\overline d(x)\subseteq \{\nu:\rho^1_m \leq \nu\in
T\}\cup (\tau \cap \sigma)\rangle_{B^c_0}$ contradiction
to (i) for $\sigma \cap \tau$).

So let $\{\tau^i:i<\alpha\}$ be maximal such that
$\overline d(\tau_i)$ are pairwise disjoint
$e^{\tau^i}_2\not= 0$, and $\tau^i\in \langle
x_\nu:\rho^1_n\leq \nu \in T\rangle_{B^c_0}$, then
$\alpha<\omega_1$, we can choose $\rho^2_n$ such that:

$\rho^1_n\leq \rho^2_n\in T$, and $[\tau\in \langle
x_\nu:\rho^2_n\leq \nu\in \tau\rangle_{B^c_0}\Rightarrow
e^\tau_2=0$.

Next we can get 

\ \ \ \ ($\beta$)  $e^\tau-1\cap e^\sigma_0=0$ (if
$\underline d(\tau)\cap\underline d(\sigma)=)$, and
$\tau,\sigma\in \langle x_\nu:\rho^2_n\leq\nu\in
T\rangle_{B^c_0})$

The proof is similar to that of $(\alpha)$, using
$\tau\cap \sigma$

As $B^c_0$ satisfies the $\aleph_1$-c.c. we can find
$\{\tau^i:i<\omega\}\subseteq \langle
x_\nu:\rho^2_n\leq\nu\in T\rangle_{B^c_0}$, such that (in
$B^c_0\, e^*_\ell\buildrel\text{def}\over =\cup_{i<\omega}
e^{\tau^i}_\ell=\cup\{ e^\tau_\ell:\tau\in \langle
x_\nu:\rho^2_n\leq\nu\in T\rangle_{B^c_0}\}$ for
$\ell=0,1$.  We can find $\rho^3_n,\rho^2_n\leq\rho^3_n\in
t$, such that $\cup_{\ell<\omega}\underline d(\tau^i)$ is
disjoint to $\{\nu:\rho^3_n\leq\nu\in T\}$.  So for every
$\tau\in \langle x_\nu\rho^3_n\leq\nu\in T\rangle_{B^c_0},
e^\tau_0\leq e^*_0$ (by the choice of $e^*_0$), and
$e^\tau_0\cap e^{\tau^i}_1=0$ for $i<\omega$ (by
($\beta$)) hence $e^\tau_0\cap e^*_1=0$, hence
 
\ \ \ \  ($\gamma$) $e^\tau_0\leq e^*_0-e^*_1$.

Similarly

\ \ \ \ ($\delta$) $e^\tau_1\leq e^*_1-e^*_0$.

Now we can prove that $e^\tau_1=e^\sigma_1$ when
$\overline d(\tau)\cap \overline d(\sigma)=0,
\tau,\sigma\in \langle x_\nu:\rho^3_n\leq\nu \in
T\rangle_{B^c_0}$, repeat the proof of $(\alpha)$
intersecting with $e^\tau_1=e^\sigma_1$ when
$\tau,\sigma\in \langle x_\nu:\rho^3_n\leq\in
T\rangle_{B^c_0}$.  So let $e_1\in B_{\alpha^*}$ be the
common value, so

$(*)\,\, h(\tau\cap y_n)\cap y_n=e^\tau_0(e_1\cap\tau)$ for
$\tau\in \langle x_\nu:\rho^3_n\leq\in T\rangle_{B^c_0}$; and
$e^\tau_0\leq y_n-e_1$,

Let $e_o=y_n=e_1$, so $y_n=e_0\cup e_1, e_0\cap e_1=0$.

So $e^\tau_0\leq e_0$ for every $\tau\in \langle
x_\nu:\rho^3_n\leq\nu \in T\rangle_{B^c_0}$

As $y_n\not\in Ex\, Ker^*(h)$, at least one of the
elements, $e_0,e_1$ is not in $Ex\, Ker(h)$>  As not
subcase III 2, for $\ell=1,2$ the homomorphism $G_\ell$
from $B_{\alpha^*}\restriction e_\ell$ to
$B_{\alpha^*}\restriction(1-e_\ell),
g_\ell(x)=h(x)-e_\ell$ (for $x=e_1$) has a finite range. 
Hence for some ideal $\Cal J$ of $B^c_0\, y_n/\Cal J$ is a
finite union of atoms and 

for every $\tau\in e\langle x_\nu:\rho^3_n\leq\nu\in
T\rangle \cap \Cal J$

for $\ell=0,1 h(\tau\cap y_n)\cap e_\ell=h(\tau\cap
e_\ell)\cap e_\ell$

hence $h(\tau \cap e_\ell)\cap e_\ell=(e^\tau_0\cup
(e^\tau_1\cap\tau))\cap e_\ell$.

So (for $\tau\in \langle x_\nu:\rho^3_n\leq\nu\in
T\rangle_{B^c_0}\cap \Cal J$:

 $h(\tau\cap e_0)\cap e_0=e^\tau_0$

  $h(\tau\cap e_1)\cap e_1=\tau\cap e_1$

If $e-1\not\in Ex\, Ker(h)$, we get contradicition to
``not subcase 3'' [use $\rho^3_n$ for $\rho^*$ there, now
for any $\rho, p^3_n\leq\rho \in T$ choose pairwise disjoint
$\tau_\ell \in \langle x_\nu:\rho\leq \in T\rangle_{B^c_0}$
for $\ell<\omega$ by the choice of $\Cal J$ for at least
one $\ell,\tau_\ell\in \Cal J$, so $\tau_\ell$ is as
requires there].  So assume $e_0\not\in Ker\, Ker^*(h)$
and get contradiction to ``not subcase III 4'' [for some
$\ell<m<\omega x_{\rho^3_n
\widehat{\,}<\ell>}-x_{\rho^3_n \widehat{\,}<n>}$) for
$\rho^*,a^*$ with $\alpha$ large enough].

So for each $n,y_n\in Ex\, Ker(h)$, by their choice) so
let $y_n=y_n^0\cup y^1_n$ (both in $B_{\alpha^*}$),
$h(y^0_n)=0, h(x)=x$ for $x\leq y^1_n, x\in
B_{\alpha^*}$.  Let $I\subseteq T$ be a countable set
such that $\underline d(y^0_n),\underline d(y^1_n)\subseteq
I$, and for $x-\in B_{\alpha^*}\, \underline
d(h(x-y_n)\cap y_n)\subseteq I$ (by ``not subcase
III 2'', for each $n$ we have only finitely many elements
of this form).

We can easily show that for every $x\in B_{\alpha^*}$, for
some $\alpha\in B^c_0$ based on $I$, $h(x)-x=a-x$, [as
$\langle y_n:n<\omega\rangle$ is a maximal antichain in
$B_{\alpha^*}$, for this it suffices to show for every
$n<\omega$ there is $a_n\in B^c_{\alpha^*}, a_n\leq y_n$
such that $(h(x)-x)\cap y_n=a_n-x$; but $(h(x)-x)\cap
y_n$ is the union of $(h(x\cap (y_n)-x\cap y_n$ which is
zeto as $(\forall_z\leq y_n)h(z)\leq z$ and of 
$(h(x-y_n)-x)\cap y_n$ which we know is bsed as wanted]. 
so $h(x)=e^x_0\cup (e^x_1\cap x)\cup e^x_2-x)$ 
where each $e^x_\ell$ is based on $I,\langle
e^x_\ell:\ell<4\rangle$ pairwise disjoint $e^*_\ell\in
B^c_0$.  As in the analysis above of $h(x\cap y_n)\cap
y_n$, (possibly with increasing $I$) applied to $x\in
B_{\alpha^*} $ with $\underline d(x)\cap I=0$, we get
$e^x_2=0, e^x_1=e_1$.  If $e_1\not\in Ex\, Ker^*(h)$ we
get contradiction to ``not subcase III 3''.  So
$1-e_1\not\in Ex\, Ker^*(h)$ and apply ``not subcase III
4''.
\endroster
So we finish the proof of 3.11; so $B_{\alpha^*}$ is
endo-rigid.

\bigskip
\proclaim{Lemma}  $B_{\alpha^*}$ is indecomposable.
\endproclaim
\bigskip\demo{Proof}  Suppose $K_0,K_1$ are disjoint
ideals of $B_{\alpha^*}$, each with no maximal members,
which generate a maximal ideal of $B_{\alpha^*}$.  For
$\ell=1,2$ let $\{d^\ell_n:\ell<\omega\}$ be a maximal
antichain $\subseteq K_\ell$ (they are countable as
$B_{\alpha^*}$ satisfies the c.c.c., and may be chosen
infinite as $K_\ell\not=\{0\}, b_{\alpha^*}$ is
atomless).  Let $K$ be the ideal $K_0\cup K_1$ generates.

Now, e.g. for some $\xi,\lambda,
\{d^\ell_n:\ell<2,n,\omega\}\subseteq B_{\xi}$.  Clearly
$a_{<\xi>}=b^0\cup b^1, b^\ell\in K_\ell$.  Now
$pr_\xi(b^\ell)\in B_{[\xi]}$ and is disjoint to each
$d^{1-\ell}_n:n<\omega\}, pr_\xi(b^\ell)$ is disjoint to
every member of $K_{1-\ell}$.  As $K-0\cup K_1$ generate
a maximal ideal, clearly $pr_\xi(b^\ell)\in K_\ell$
[otherwise $pr_\xi(b^\ell)=1-c^1\cup c^2$, for some
$c^1\in K_1, c^2\in K_2$, and then $c^{1-\ell}$ is
necessarily a maximal member of $K_{1-\ell}$, so
$K_{1-\ell}$ is principal contradiction].  So
$pr_\xi(B^0)\cup pr_\xi(b^2)<1$ but
$1=pr_\xi(a_{<\xi>}=\cup^2_{\ell=0} pr_\xi(b^\ell)$
contrdiction.\enddemo
\bigskip
\proclaim{3.13 Theorem}  In 3.1 we can get
$2{\lambda^{\aleph_0}}$ such Boolean Algebras  such that
any homomorphism from one to the other has a finite range.
\endproclaim
\bigskip
\demo{Proof}  Left to the reader (see [Sh 4.3]).
\enddemo
\vfill
\pagebreak

\enddocument